\documentclass[a4paper,11pt]{article}
\usepackage{amsmath,amscd,amsfonts,amssymb,epsf,latexsym}
\usepackage{colordvi,color}     %%  ADDED
\makeatletter

\long\def\@makefnt#1{\parindent 1em\noindent
            \hb@xt@1.8em{\hss\@textsuperscript{}}#1}
\long\def\@ftntext#1{\insert\footins{%
    \reset@font\footnotesize
    \interlinepenalty\interfootnotelinepenalty
    \splittopskip\footnotesep
    \splitmaxdepth \dp\strutbox \floatingpenalty \@MM
    \hsize\columnwidth \@parboxrestore
    \color@begingroup
      \@makefnt{%
        \rule\z@\footnotesep\ignorespaces#1\@finalstrut\strutbox}%
    \color@endgroup}}%
\def\subjclass#1{%
  \@ftntext{2000 {\itshape Mathematics Subject Classification.}\enspace #1.}}
\def\keywords#1{%
  \@ftntext{{\itshape Key words and phrases.}\enspace #1.}}
\makeatother

\def\A{{\mathbb A}}
\def\B{{\mathbb B}}
\def\N{{\mathbb N}}
\def\C{{\mathbb C}}
\def\D{{\mathbb D}}

\def\X{{\mathbb X}}
\def\Y{{\mathbb Y}}

\def\AB{ {\mathbb A}\moins {\mathbb B}}

\def\moins{\raise 1pt\hbox{{$\scriptstyle -$}}}
\def\plus{\raise 1pt\hbox{{$\scriptstyle +$}} }

\def\phi{\varphi}

\newtheorem{theorem}{Theorem}
\newtheorem{proposition}[theorem]{Proposition}
\newtheorem{lemma}[theorem]{Lemma}
\newtheorem{corollary}[theorem]{Corollary}
\newtheorem{remark}[theorem]{Remark}
\newtheorem{definition}[theorem]{Definition}

\def\proof{\noindent{\bf Proof.\ }}

\def\qed{~\hbox{$\Box$}}

\def\Aut{\mathop{\rm Aut}}
\def\Diff{\mathop{\rm Diff}}
\def\codim{\mathop{\rm codim}}

\def\dim{\mathop{\rm dim}}

\begin{document}

\title{\bf Thom polynomials and Schur functions: the singularities $A_3(-)$}

\author{Alain Lascoux\thanks{Research supported by the ANR project MARS (BLAN06-2 134516).}\\
\small IGM, Universit\'e de Paris-Est\\
\small 77454 Marne-la-Vall\'ee CEDEX 2\\
\small Alain.Lascoux@univ-mlv.fr \and
Piotr Pragacz\thanks{Research partially supported by the University of Kyoto during the
author's stay at RIMS.}\\
\small Institute of Mathematics of Polish Academy of Sciences\\
\small \'Sniadeckich 8, 00-956 Warszawa, Poland\\
\small P.Pragacz@impan.gov.pl}

\subjclass{05E05, 14N10, 57R45}

\keywords{Thom polynomials, singularities, global singularity
theory, classes of degeneracy loci, Schur functions, resultants, Pascal staircases}

\date{}

\maketitle

\centerline{\it To the memory of Stanis\l aw Balcerzyk}

\begin{abstract}
Combining the ``method of restriction equations'' of Rim\'anyi et al. with the techniques 
of symmetric functions, we establish the Schur function expansions 
of the Thom polynomials for the Morin singularities 
$A_3: ({\bf C}^{\bullet},0)\to ({\bf C}^{\bullet + k},0)$
for any nonnegative integer $k$.
\end{abstract}

\section{Introduction}

The global behavior of singularities of maps is governed by their {\it Thom
polynomials} (see \cite{T}, \cite{Kl}, \cite{AVGL}, \cite{Ka},
\cite{Rim2}). Knowing the Thom polynomial of a singularity $\eta$,
denoted ${\cal T}^{\eta}$, one can compute the cohomology class
represented by the $\eta$-points of a map. In particular, if $f:X\to Y$
is a general map of complex analytic manifolds, where $X$ is compact and $\dim(X)$
equals the codimension of the singularity $\eta$, then the degree $\int_{X} {\cal T}^{\eta}$ 
evaluates the number of points of $X$ at which $f$ has the singularity $\eta$.

In the present paper, following the ``method of restriction equations'' 
from a series of papers by Rim\'anyi et al. \cite{RS}, \cite{Rim2}, \cite{FR}, \cite{BFR},
we study the Thom polynomials for the singularities
$A_3$ associated with maps $({\bf C}^{\bullet},0)
\to ({\bf C}^{\bullet+k},0)$  with parameter $k\ge 0$.
We give the Schur function expansions of these Thom polynomials.
This is the content of our main Theorem \ref{TA3} and its proof 
in Section \ref{result}.

The way of obtaining the Thom polynomial is through the solution
of a system of linear equations (see Theorem \ref{TEq}). This is fine when 
we want to find one concrete Thom polynomial, say, for a fixed $k$. However, if we want
to find the Thom polynomials for a series of singularities, associated
with maps $({\bf C}^{\bullet},0) \to ({\bf C}^{\bullet+k},0)$
with $k$ as a parameter, we have to solve {\it simultaneously} a countable
family of systems of linear equations. This cannot be done by computer,
and must be done conceptually.

Thom polynomials are symmetric functions in the universal Chern roots. Instead of
giving their expressions in terms of these variables, we use
{\it Schur function expansions}. This puts a more transparent structure 
on computations of Thom polynomials (see \cite{P23}, and also \cite{FK} 
for some second order Thom-Boardman singularities). In particular, in the 
Schur basis one can see some {\it recurrences} which are difficult or even 
impossible to notice in other bases (see \cite{P4}).

Another feature of using the Schur function expansions for Thom
polynomials is that all the coefficients are {\it nonnegative}.
This has been recently proved by A.~Weber and the second author in \cite{PW}. 

To be more precise, we use here (the specializations of)
{\it supersymmetric} Schur functions also called
``Schur functions in difference of alphabets'' together with their
three basic properties: {\it vanishing}, {\it cancellation} and {\it
factorization}, (see \cite{S}, \cite{BR}, \cite{P2}, \cite {PT},
\cite{M}, and \cite{L}). These functions contain resultants among themselves. 
They play a fundamental role in the study of 
${\cal P}$-{\it ideals of singularities $\Sigma^i$}
(see \cite[end of Sect.~2 and Theorem 11]{P4} and Proposition \ref{sigmai}) which 
is based on the enumerative geometry of degeneracy loci of \cite{P}.

Since the singularity $A_3$ is in the closure of the orbit of the singularity 
$\Sigma^1$, we have by Proposition \ref{sigmai} that all partitions in the Schur 
expansion of ${\cal T}^{A_3}$ (any $k$) contain the single row-partition $(k+1)$.

In \cite{P3} the decomposition of the Thom polynomial of the singularity $A_i$ into
$h$-{\it parts} was defined (see also the end of Section \ref{Schur}). In particular,
the $1$-part of the Thom polynomial of the Morin singularity $A_i$ (any $i, k$)
was computed. In the present paper, we work out the case of the singularities $A_3$
(any $k$), and we find the $2$-part of this Thom polynomial (the $h$-parts, 
where $h\ge 3$, are equal to zero for these singularities).

\smallskip

In our calculations, we use extensively the functorial $\lambda$-ring
approach to symmetric functions from \cite{L} (e.g. we shall need 
to handle symmetric functions in $2x_1,2x_2,x_1+x_2$\footnote{Strictly speaking: 
symmetric functions in $\fbox{$2x_1$}$, $\fbox{$2x_2$}$, $\fbox{$x_1+x_2$}$ after
simplification, see Section \ref{Schur}.}
at the same time as symmetric functions in $x_1,x_2$).

\smallskip

The main results of the present paper were announced in \cite{P23}. 

\smallskip

B\'erczi, Feh\'er and R. Rim\'anyi gave without proof in 
\cite{BFR} an expression for this Thom polynomial, but in 
terms of the monomial basis in Chern classes. We prompt
the authors of \cite{BFR} to publish their proof.

\section{Reminder on Thom polynomials}

Our main reference for this section is \cite{Rim2}.
We start with recalling what we shall mean by a ``singularity''.
Let $k\ge 0$ be a fixed integer. By a {\it singularity}
we shall mean an equivalence class of stable germs $({\bf C}^{\bullet},0)
\to ({\bf C}^{\bullet+k},0)$, where $\bullet\in {\bf N}$, under the
equivalence
generated by right-left equivalence (i.e. analytic reparametrizations
of the source and target) and suspension.

We recall\footnote{This statement is usually called the Thom-Damon
theorem \cite{T}, \cite{D}.}
that the {\it Thom polynomial} ${\cal T}^{\eta}$ of a singularity
$\eta$ is a polynomial in the formal variables $c_1, c_2, \dots$ which
after the substitution of $c_i$ to
\begin{equation}\label{ci}
c_i(f^*TY-TX)=[c(f^*TY)/c(TX)]_i\,,
\end{equation}
for a general map $f:X \to Y$ between complex analytic manifolds,
evaluates the Poincar\'e dual of $[V^{\eta}(f)]$, where $V^{\eta}(f)$
is the cycle carried by the closure of the set
\begin{equation}
\{x\in X : \hbox{the singularity of} \ f \ \hbox{at} \ x \
\hbox{is} \ \eta \}\,.
\end{equation}
By {\it codimension of a singularity} $\eta$, $\codim(\eta)$,
we shall mean $\codim(V^{\eta}(f),X)$ for such an $f$. The concept of
the polynomial ${\cal T}^{\eta}$ comes from Thom's fundamental paper
\cite{T}.
For a detailed discussion of the {\it existence} of Thom polynomials,
see, e.g., \cite{AVGL}. Thom polynomials associated with group actions
were studied by Kazarian in \cite{Ka} and \cite{Ka2}.

According to Mather's classification, singularities are in one-to-one
correspondence with finite dimensional ${\bf C}$-algebras.
We shall use the following notation:

\begin{itemize}
\def\labelitemi{--}
\item $A_i$ \ (of Thom-Boardman type $\Sigma^{1_i}$) will stand for the stable
germs with local algebra ${\bf C}[[x]]/(x^{i+1})$, $i\ge 0$;

\item $III_{2,2}$ \ (of Thom-Boardman type $\Sigma^2$) for stable germs
with local algebra ${\bf C}[[x,y]]/(xy, x^2, y^2)$ (here $k\ge 1$).
\end{itemize}

In the present article, the computations of Thom polynomials
shall use the method which stems from a sequence
of papers by Rim\'anyi et al. \cite{RS}, \cite{Rim2},
\cite {FR}, \cite{BFR}.
We sketch briefly this approach, refering the interested reader for more
details to these papers.

Let $k\ge 0$ be a fixed integer, and let $\eta: ({\bf C}^{\bullet},0) \to
({\bf C}^{\bullet+k},0)$ be a stable singularity with a prototype
$\kappa: ({\bf C}^n,0) \to ({\bf C}^{n+k},0)$. The {\it maximal compact
subgroup of the right-left symmetry group}
\begin{equation}
\Aut \kappa = \{(\phi,\psi) \in \Diff({\bf C}^n,0) \times
\Diff({\bf C}^{n+k},0) : \psi \circ \kappa \circ \phi^{-1} = \kappa \}
\end{equation}
of $\kappa$ will be denoted by $G_\eta$.
Even if $\Aut \kappa$ is much too large to be a finite dimensional
Lie group, the concept of its maximal compact subgroup (up to conjugacy)
can be defined in a sensible way.
In fact, $G_\eta$ can be chosen so that the images of its projections
to the factors $\Diff({\bf C}^n,0)$ and $\Diff({\bf C}^{n+k},0)$ are
linear. Its representations via the projections on the source ${\bf C}^n$
and the target ${\bf C}^{n+k}$ will be denoted by $\lambda_1(\eta)$ and
$\lambda_2(\eta)$.
The vector bundles associated with the universal principal
$G_\eta$-bundle $EG_\eta \to BG_\eta$ using the representations
$\lambda_1(\eta)$ and $\lambda_2(\eta)$ will be called
$E_{\eta}'$ and $E_{\eta}$. The {\it total Chern
class of the singularity} $\eta$ is defined
in $H^{*}(BG_\eta, {\bf Z})$ by
\begin{equation}
c(\eta):=\frac{c(E_{\eta})}{c(E_{\eta}')}\,.
\end{equation}
The {\it Euler class} of $\eta$ is defined in
$H^{2\codim(\eta)}(BG_\eta, {\bf Z})$ by
\begin{equation}
e(\eta):=e(E_{\eta}')\,.
\end{equation}

In the following theorem, we collect information from \cite{Rim2},
Theorem 2.4 and \cite{FR}, Theorem 3.5, needed for the calculations
in the present paper.

\begin{theorem}\label{TEq} Suppose, for a singularity $\eta$, that
the Euler classes of all singularities of smaller codimension than
$\codim(\eta)$, are not zero-divisors \footnote{This is the so-called
``Euler condition'' ({\it loc.cit.}). It holds for $A_3$.}. Then we have

\begin{enumerate}
\item  if $\xi\ne \eta$ and  $\codim(\xi)\le \codim(\eta)$, then
 ${\cal T}^{\eta}(c(\xi))=0${\rm;}

\item  ${\cal T}^{\eta}(c(\eta))=e(\eta)$.

\end{enumerate}
This system of equations $($taken for all such $\xi$'s$)$ determines the
Thom polynomial ${\cal T}^{\eta}$ in a unique way%
\footnote{To make it precise, we need one more condition
that the number of singularities (=contact orbits) of smaller codimension
is finite: we may assume that $\eta$ is a {\it simple} singularity type,
i.e., there is no moduli adjacent to $\eta$.}.
\end{theorem}

To use this method of determining the Thom polynomials for
singularities, one needs their classification, see, e.g., \cite{dPW}.

We record the following lemma (see \cite{Rim2} and \cite{BFR}). 

\begin{lemma}\label{ce} \ (i) For the singularity of type 
$A_i$: $({\bf C}^{\bullet},0) \to ({\bf C}^{\bullet+k},0)$, we have 
$G_{A_i}=U(1)\times U(k)$. Moreover, denoting by $x$ and $y_1,\ldots,y_k$
the Chern roots of the tautological vector bundles on $BU(1)$ and $BU(k)$,
we have
\begin{equation}\label{cAi}
c(A_i)=\frac{1+(i+1)x}{1+x}\prod_{j=1}^k (1+y_j)
\end{equation}
and
\begin{equation}\label{eA}
e(A_3)= 6 \ x^3\ \prod_{j=1}^k (y_j-3x)(y_j-2x)(y_j-x)\,.
\end{equation}

\noindent
(ii) For the singularity $III_{2,2}: ({\bf C}^{\bullet},0) \to ({\bf C}^{\bullet+k},0)$,
where $k>0$, we have $G_\eta=U(2)\times U(k\moins 1)$. Moreover, 
denoting by $x_1,x_2$ (resp. $y_1,\ldots, y_{k-1}$) the Chern roots of the tautological vector
bundle on $BU(2)$ (resp. $BU(k-1)$), we have

\begin{equation}\label{eIII}
c(III_{2,2})=\frac{(1\plus 2x_1)(1\plus 2x_2)(1\plus x_1\plus x_2)}
{(1\plus x_1)(1\plus x_2)}\prod_{j=1}^{k-1}(1+y_j)\,.
\end{equation}
\end{lemma}

\section{Reminder on Schur functions}\label{Schur}

In this section, we collect needed notions related to symmetric
functions. We adopt a functorial $\lambda$-ring point of view of \cite{L}.

\smallskip

For $m\in {\bf N}$, by an alphabet $\A$ we shall mean a finite set of
indeterminates $\A=\{a_1,\dots,a_m\}$.

We shall often identify an alphabet
$\A=\{a_1,\dots,a_m\}$ with the sum $a_1+\cdots +a_m$.

\begin{definition}\label{cf}
Given two alphabets $\A$, $\B$, the {\it complete functions} $S_i(\AB)$
are defined by the generating series (with $z$ an extra variable):

\begin{equation}
\sum S_i(\AB) z^i =\prod_{b\in \B} (1\moins bz)/\prod_{a\in \A}
(1\moins az)\,.
\end{equation}
\end{definition}

\begin{definition}\label{sf}
Given a partition\footnote{We identify partitions with their
Young diagrams, as is customary.} $I=(0\le i_1\le i_2\le \ldots\le i_s)\in
{\bf N}^s$, and alphabets $\A$ and $\B$, the {\it Schur function}
$S_I(\A \moins \B)$ is \begin{equation}\label{schur}
S_I(\A \moins \B):= \bigl|
     S_{i_p+p-q}(\A \moins \B) \bigr|_{1\leq p,q \le s}  \, .
\end{equation}
\end{definition}

These functions are often called {\it supersymmetric Schur functions}
or {\it Schur functions in difference of alphabets}. Their properties
were studied, among others, in \cite{BR}, \cite{P2},
\cite {PT}, \cite{M}, and \cite{L}.

We have the following {\it cancellation property}:
\begin{equation}
S_I((\A + \C) - (\B + \C))=S_I(\A-\B)\,.
\end{equation}

We identify partitions with their Young diagrams, as is customary.

We record the following property ({\it loc.cit.}):
\begin{equation}
S_I(\AB)= (-1)^{|I|}S_J(\B \moins \A)=S_J(\B^* \moins \A^*)\,,
\end{equation}
where $J$ is the conjugate partition of $I$ (i.e. the consecutive
rows of the diagram of~$J$ are the transposed columns of the
diagram of $I$), and
$\A^*$ denotes the alphabet $\{-a_1,-a_2,\dots \}$.

%Fix two positive integers $m$ and $n$.
%We shall say that a partition $I=(0<i_1\le i_2\le \cdots \le i_s)$
%{\it is contained in} the $(m,n)$-hook if either $s\le m$, or $s> m$
%and $i_{s-m}\le n$.
%Pictorially, this means that the Young diagram of $I$ is contained
%in the hook:

%\smallskip

%$$
%\unitlength=2mm
%\begin{picture}(18,14)
%\put(0,0){\line(0,1){14}}
%\put(0,0){\line(1,0){18}}
%\put(9,5){\line(1,0){9}}
%\put(9,5){\line(0,1){9}}
%\put(4,10){\vector(1,0)5}
%\put(4,10){\vector(-1,0)4}
%\put(13,2){\vector(0,1)3}
%\put(13,2){\vector(0,-1)2}
%\put(4.5,10.3){\hbox to0pt{\hss$n$\hss}}
%\put(13.3,2.3){\hbox{$m$}}
%\end{picture}
%$$

%\smallskip

%We record the following {\it vanishing property}.
%Given alphabets $\A$ and $\B$ of cardinalities $m$ and $n$, if
%a partition $I$ is not contained in the $(m,n)$-hook,
%then ({\it loc.cit.}):
%\begin{equation}\label{van}
%S_I(\A-\B)=0\,.
%\end{equation}

\smallskip

In the present paper, by a {\it symmetric function} we shall mean
a ${\bf Z}$-linear combination of the operators $S_I$.

Instead of introducing, in the argument of a symmetric function,
formal variables which will be specialized, we write
$\fbox{$r$}$ for a variable which will be specialized to $r$
($r$ can be $2x_1$, $x_1+x_2$,\ldots).
For example,
$$
S_2(x_1+x_2)=x_1^2+ x_1x_2+ x_2^2 \ \ \ \hbox{but} \ \ \
S_2\bigl(\fbox{$x_1\plus x_2$}\bigr)=
(x_1+x_2)^2= x_1^2+ 2x_1x_2+ x_2^2\,.
$$

\begin{definition}
Given two alphabets $\A,\B$, we define their {\it resultant}:
\begin{equation}\label{res}
R(\A,\B):=\prod_{a\in \A,\, b\in \B}(a\moins b)\,.
\end{equation}
\end{definition}

For example, we have the following identity:
\begin{equation}\label{fid}
-6x^3\prod_{j=1}^k (3x\moins y_j)(2x\moins y_j)(x\moins y_j)=
R\bigl(x\plus \fbox{$2x$}\plus \fbox{$3x$},\Y \plus \fbox{$4x$}\bigr)\,,
\end{equation}
where $\Y=\{y_1,\ldots,y_k\}$.

%We have (see \cite{L})
%\begin{equation}\label{ER}
%R(\A_m,\B_n)= S_{(n^m)}(\AB)=\sum_I S_I(\A) S_{(n^m)/I}(-\B)\,,
%\end{equation}
%where the sum is over all partitions $I\subset (n^m)$.

\smallskip

We record the following {\it factorization property} (\cite[Proposition 1.4.3]{L}).
Suppose that cardinality of $B$ is $n$. Then for partitions $I=(i_1,\dots,i_m)$ and 
$J=(j_1,\dots, j_s)$, we have
\begin{equation}\label{Fact}
S_{(j_1,\dots,j_s,i_1+n,\dots,i_m+n)}(\A-\B)
=S_I(\A) \ R(\A,\B) \ S_J(-\B)\,.
\end{equation}

\smallskip

In the present paper, it will be more handy to use, instead of $k$, a shifted 
parameter
\begin{equation}
r:=k+1\,.
\end{equation}
Sometimes, we shall write $\eta(r)$ for the singularity
$\eta: ({\bf C}^{\bullet},0) \to ({\bf C}^{\bullet + r-1},0)$,
and denote the Thom polynomial of $\eta(r)$ by ${\cal T}^{\eta}_r$
-- to emphasize the dependence of both items on $r$.

\smallskip

Let $f:X\to Y$ be a map of complex analytic manifolds, where $\dim(X)=m$
and $\dim(Y)=n$. Given a partition $I$, we define 
$$
S_I(T^*X-f^*(T^*Y))
$$
to be the effect of the following specialization of $S_I(\A-\B)$: we set
the indeterminates of $\A$ to the Chern roots of $T^*X$, and the indeterminates
of $\B$ to the Chern roots of $f^*(T^*Y)$. 

Similarly to \cite{P23}, \cite{P4}, and \cite{P3}, we shall write the Poincar\'e
dual of $[V^{\eta}(f)]$, for a singularity $\eta$ and a general map $f: X\to Y$,
in the form 
$$
\sum_I \alpha_I S_I(T^*X-f^*(T^*Y))
$$
with integer coefficients $\alpha_I$. Accordingly, we shall write
\begin{equation}
{\cal T}^{\eta}=\sum_I \alpha_I S_I\,,
\end{equation}
where $S_I$ is identified with $S_I(\A-\B)$ for the universal Chern roots
$\A$ and $\B$. 

Note that in this notation, the Thom polynomial of the singularity
$A_1(r)$ for $r\ge 1$, is: ${\cal T}^{A_1}_r=S_r$.
Another example is the Thom polynomial of $A_2(1)$. In \cite{Rim2}, it is written as
$c_1^2+c_2$, whereas in the present notation it is written as $S_{11}+2S_2$.

\medskip

The arguments of the proof of \cite[Theorem 11]{P4} give the following result\footnote{This 
justifies the remark in \cite{P3} p.~166, lines 28--31.}.

\begin{proposition}\label{sigmai} \ Suppose that a singularity $\eta$ is in the closure 
of the orbit of the singularity $\Sigma^j$. Then all summands in the Schur function expansion 
of ${\cal T}^{\eta}_r$ are indexed by partitions containing the rectangle partition $(r+j-1)^j$.
\end{proposition}

\medskip

Recall (from \cite{P3}) that the $h$-{\it part} of ${\cal T}^{A_i}_r$
is the sum of all Schur functions appearing nontrivially
in ${\cal T}^{A_i}_r$ (multiplied by their coefficients) such that the
corresponding partitions satisfy the following condition:
$I$ contains the rectangle partition $\bigl((r\plus h\moins 1)^h\bigr)$,
but it does not contain the larger Young diagram $\bigl((r\plus h)^{h+1}\bigr)$.
The polynomial ${\cal T}^{A_i}_r$ is a sum of its $h$-parts, $h=1,2,\dots$.

\medskip

In one instance (the proof of Proposition \ref{gr}), we shall also use 
{\it multi-Schur functions}. 
For their definition and properties, we refer the reader to \cite{L}.

\section{Main result and its proof}\label{result}

Since the singularities $\ne A_3$, whose codimension is $\le \codim(A_3)$
are: $A_0$, $A_1$, $A_2$ and, for $r\ge 2$, $III_{2,2}$ (see \cite{dPW}),
Theorem \ref{TEq} yields the following equations (in $T$),
characterizing the Thom polynomial ${\cal T}^{A_3}_r$:
\begin{equation}\label{EqA3}
T(-\B_{r-1})=T(x-\B_{r-1}-\fbox{$2x$})=T(x\moins \B_{r-1}-\fbox{$3x$})=0\,,
\end{equation}
\begin{equation}\label{EqA3'}
T(x- \B_{r-1}-\fbox{$4x$})=
R(x+ \fbox{$2x$}+ \fbox{$3x$}, \B_{r-1}+ \fbox{$4x$}\, )
\end{equation}
\begin{equation}
T(x_1+x_2-\D-\B_{r-2})=0\,.
\end{equation}
Here,
$$
\D=\fbox{$2x_1$}+\fbox{$2x_2$}+\fbox{$x_1+x_2$}\,.
$$

We assume that $x$, $x_1$, $x_2$,
and $b_1,\ldots, b_n$ are variables.
Note that these variables, in the following, will be specialized to the Chern roots 
of the {\it cotangent} bundles. 

\smallskip

By \cite {P3}, we know that ${\cal T}^{A_3}_r$ must contain (as its $1$-part)
the following combination of Schur functions, denoted by $F_r^{(3)}$ in \cite{P3}:

\begin{equation}
F_r:=\sum_{j_1\le j_2 \le r} S_{j_1,j_2}(\fbox{$2$}+\fbox{$3$})
S_{r-j_2,r-j_1,r+j_1+j_2}\,.
\end{equation}

\smallskip

By \cite[Corollary 11]{P3}, Eqs. (\ref{EqA3}) and (\ref{EqA3'}) are satisfied by
the function $F_r$. For $r=1$, this means that
\begin{equation}
F_1=S_{111}+5S_{12}+6S_3
\end{equation}
is the Thom polynomial for $A_3(1)$.

\smallskip

However, for $r\ge 2$, $F_r$ does not satisfy the last vanishing,
imposed by $III_{2,2}$.
In the following we shall modify $F_r$ in order to obtain
the Thom polynomial for $A_3$.
In fact, our goal is to give an expresion for the Thom polynomial for
$A_3$ (any $r$) as a ${\bf Z}$-linear combination of Schur functions.
For $r=2$, the Thom polynomial is
\begin{equation}
S_{222}+5S_{123}+6S_{114}+19S_{24}+30S_{15}+36S_6+5S_{33}\,,
\end{equation}
and it differs from its $1$-part $F_2$ by $5S_{33}$ which is
the ``correction'' $2$-part in this case (see \cite{P3}).

\smallskip

Define integers $e_{i,j}$, for $i\ge 2$ and $j\ge 0$ in the following way.
First, $e_{20}, e_{30}, e_{40},\ldots$ are the coefficients
$5, 24, 89,\ldots$ in the Taylor expansion of
$$
\aligned
\frac{5-6z}{(1-z)(1-2z)(1-3z)}\\
=5+24z+&89z^2+300z^3+965z^4+3024z^5+9329z^6+\ldots\,.
\endaligned
$$
Moreover, we set $e_{2,j}=e_{3,j}=0$ for $j\ge 1$, $e_{4,j}=e_{5,j}=0$
for $j\ge 2$, $e_{6,j}=e_{7,j}=0$ for $j\ge 3$ etc.
To define the remaining $e_{i,j}$'s, we use the recursive formula
\begin{equation}\label{re}
e_{i+1,j}= e_{i,j-1} + e_{i,j}\,.
\end{equation}
We obtain the following matrix $[e_{i,j}]_{i\ge 2, j\ge 0}$~:
$$
\begin{array}{cccccc}
e_{20} & 0 & 0 & 0 & 0 & \ldots \\
e_{30} & 0 & 0 & 0 & 0 & \ldots \\
e_{40} & e_{41} & 0 & 0 & 0 & \ldots \\
e_{50} & e_{51} & 0 & 0 & 0 & \ldots \\
e_{60} & e_{61} & e_{62} & 0 & 0 & \ldots \\
e_{70} & e_{71} & e_{72} & 0 & 0 & \ldots \\
e_{80} & e_{81} & e_{82} & e_{83} & 0 & \ldots \\
\vdots & \vdots & \vdots & \vdots & \vdots &
\end{array} \ \ \ \ \ \
= \ \ \ \ \ \
\begin{array}{cccccc}
5 & 0 & 0 & 0 & 0 & \ \ldots \\
24 & 0 & 0 & 0 & 0 & \ldots \\
89 & 24 & 0 & 0 & 0 & \ldots \\
300 & 113 & 0 & 0 & 0 & \ldots \\
965 & 413 & 113 & 0 & 0 & \ldots \\
3024 & 1378 & 526 & 0 & 0 & \ldots \\
9329 & 4402 & 1904 & 526 & 0 & \ldots \\
\vdots & \vdots & \vdots & \vdots & \vdots &
\end{array}
$$

\begin{remark} \rm 
Note that arguing similarly as in the proof of Proposition 19 in \cite{P4}, we
get the following closed formula for $e_{i,j}$.
For $i\ge 2$ and $j\ge 0$, we have
$$
\aligned
e_{i,j}=\frac{1}{2^{j\plus 1}}\Bigl[(3^{i\plus 1}\moins 3^{2(j\plus 1)})
&\moins (2^{i\plus j\plus 2}-2^{3(j\plus 1)})\\
\moins \sum_{s=1}^j 2^s\bigl(3^{2(j\moins s\plus 1)}
&\moins 2^{3(j\moins s\plus 1)}\bigr)
\Bigl({i\moins 2j\moins 2s\plus 1 \choose s}
\moins {2s\moins 2 \choose s} \Bigr)\Bigr]\,.
\endaligned
$$
For example, we have
$$
e_{i,2}=\frac{1}{2^3}[(3^{i\plus 1}\moins 3^6)\moins (2^{i\plus 4}\moins 2^9)
\moins 2(3^4\moins 2^6)(i\moins 5)
\moins 2^2\Bigl({i\moins 3 \choose 2}\moins 1\Bigr)\Bigr]\,.
$$
\end{remark}

\medskip

Consider the following matrix whose elements are two row partitions
(the symbol ``$\emptyset$" denotes the empty partition):
$$
\begin{array}{cccccc}
33 & \emptyset & \emptyset & \emptyset & \emptyset & \ldots \\
45 & \emptyset & \emptyset & \emptyset & \emptyset & \ldots \\
57 & 66 & \emptyset & \emptyset & \emptyset & \ldots \\
69 & 78 & \emptyset & \emptyset & \emptyset & \ldots \\
7,11 & 8,10 & 9,9 & \emptyset & \emptyset & \ldots \\
8,13 & 9,12 & 10,11 & \emptyset & \emptyset & \ldots \\
9,15 & 10,14 & 11,13 & 12,12 & \emptyset & \ldots \\
\vdots & \vdots & \vdots & \vdots & \vdots &
\end{array}
$$
\noindent
We use for this matrix the same ``matrix coordinates'' as for the
previous one.
Denote by $I(i,j)$ the partition occupying the $(i,j)$th place
in this matrix. So, e.g., $I(i,0)=(i+1,2i-1)$ for $i\ge 2$. 

\medskip

For $r\ge 2$, we set
\begin{equation}
{\overline H_r}:=\sum_{j\ge 0} e_{r,j} \ S_{I(r,j)}\,.
\end{equation}
We have

$$\begin{array}{lllll}
{\overline H_2}&=5S_{33}\\
{\overline H_3}&=24S_{45}\\
{\overline H_4}&=89S_{57} &+24S_{66}\\
{\overline H_5}&=300S_{69} &+113S_{78}\\
{\overline H_6}&=965S_{7,11} &+413S_{8,10} &+113S_{99}\\
{\overline H_7}&=3024S_{8,13} &+1378S_{9,12} &+526S_{10,11}\\
{\overline H_8}&=9329S_{9,15} &+4402S_{10,14} &+1904S_{11,13} &+526S_{12,12}\,.
\end{array} $$

Denote now by $\Phi$ the linear endomorphism on the free ${\bf Z}$-module
spanned by Schur functions indexed by partitions of length $\le 3$,
that sends a Schur function $S_{i_1,i_2,i_3}$ to $S_{i_1+1,i_2+1,i_3+1}$.
We define
\begin{equation}\label{recH}
H_r:={\overline H_r}+\Phi(H_{r-1})\,,
\end{equation}
or equivalently, by iteration
\begin{equation}
H_r={\overline H_r}+\Phi({\overline H_{r-1}})+\Phi^2({\overline H_{r-2}})
+\cdots+\Phi^{r-2}({\overline H_2})\,.
\end{equation}
We have the following values of
$H_2, H_3=\Phi(H_2)+\overline{H_3},\ldots, H_7=\Phi(H_6)+\overline{H_7}$~:
$$
\aligned
H_2&=5S_{33}\\
H_3&=5S_{144}\plus24S_{45}\\
H_4&=5S_{255}\plus24S_{156}\plus24S_{66}\plus89S_{57}\\
H_5&=5S_{366}\plus24S_{267}\plus24S_{177}\plus89S_{168}\plus113S_{78}
\plus300S_{69},\\
H_6&=5S_{477}\plus24S_{378}\plus24S_{288}\plus89S_{279}\plus113S_{189}
\plus300S_{1,7,10}\plus113S_{99}\plus413S_{8,10}\\
&\plus965S_{7,11}\\
H_7&=5S_{588}\plus24S_{489}\plus24S_{399}\plus89S_{3,8,10}\plus113S_{2,9,10}
\plus300S_{2,8,11}\plus113S_{1,10,10}\\
&\plus413S_{1,9,11}\plus 965S_{1,8,12}\plus526S_{10,11}\plus1378S_{9,12}\plus3024S_{8,13}\,.
\endaligned
$$

Alternatively, 

\begin{equation}\label{hr}
H_r = \sum_{i=0}^{r-2} \ \ \sum_{\{j\ge 0: \ i+2j \le r-2\}}
 \ e_{r-i,j} \ S_{i,r+j+1,2r-i-j-1}\,.
\end{equation}

\medskip

We now state the main result of this paper.

\begin{theorem}\label{TA3}
For $r\ge 1$, the Thom polynomial of $A_3(r)$ is equal to $F_r+H_r$.
\end{theorem}
In other words, the function $H_r$ is the $2$-part of ${\cal T}^{A_3}_r$, 
and its $h$-parts are zero for $h\ge 3$.

\medskip

In the proof of the theorem, we shall need several properties of
the functions $H_r$ and $F_r$.

\medskip

The next result says that the addition of $H_r$ to $F_r$
is ``irrelevant'' for what concerns the conditions (\ref{EqA3}) and (\ref{EqA3'}) 
imposed by the singularities $A_i$, $i=0,1,2,3$.

\begin{lemma}\label{Lv} \ The function $H_r$ satisfies
Eqs.~(\ref{EqA3}), and the equation
\begin{equation}
H_r(x-\B_{r-1}-\fbox{$4x$})=0\,.
\end{equation}
\end{lemma}
\proof
According to (\ref{Fact}), each Schur function of index $(i_1,i_2,i_3)$ with 
$i_2,i_3\geq r\plus 1$ vanishes when evaluated in $x-\B_{r-1}-y$, $y$ any indeterminate.
Therefore $H_r$ satisfies the required nullities, which correspond to
taking $y=0,x, \fbox{$2x$}, \fbox{$3x$}$ or $\fbox{$4x$}$.
\qed

\smallskip

Thanks to the lemma, in order to prove the theorem, it suffices to show the equality
\begin{equation}\label{frhr}
(F_r+H_r)(x_1+x_2-\D-\B_{r-2})=0\,,
\end{equation}
which is equivalent to the vanishing of ${\cal T}^{A_3}_r$ at the Chern class $c(III_{2,2}(r))$.

\medskip

Set $\X_2=(x_1,x_2)$. Due to (\ref{Fact}), each Schur function occuring in the expansion of $H_r$ 
is such that

\begin{equation*}
S_{c,r+1+a,r+1+b}(\X_2\moins \D\moins \B_{r-2})
= R(\X_2,\D\plus \B_{r-2})\cdot S_c(\moins \D\moins \B_{r-2})\cdot S_{a,b}(\X_2)\,,
\end{equation*}
We set
\begin{equation}
V_r(\X_2;\B_{r-2})= \frac{H_r(\X_2-\D-\B_{r-2})}{R(\X_2,\D+\B_{r-2})}\,,
\end{equation}
so that
\begin{equation}\label{v}
V_r(\X_2;\B_{r-2})=\sum_{i=0}^{r-2} \sum_{\{j\ge 0: \ i+2j \le r-2\}}
e_{r-i,j} \ S_i(-\D-\B_{r-2}) \ S_{j,r-i-j-2}(\X_2)\,.
\end{equation}
We have the following
recursive relation which follows from the observation that the coefficient
of $b_{r-2}$ in $V_r(\X_2;\B_{r-2})$ is equal to
$-V_{r-1}(\X_2;\B_{r-3})$.

\begin{lemma}\label{Lvr} For $r\ge 2$, we have
\begin{equation}
V_r(\X_2; \B_{r-2})= \sum_{i=0}^{r-2} \ V_{r-i}(\X_2; 0)
\ S_i(-\B_{r-2})\,.
\end{equation}
\end{lemma}
Thus it is sufficient to compute $V_r(\X_2;0)$.

\begin{proposition}\label{Pv} \ For $r\ge 2$, we have
\begin{equation}\label{Vr0}
V_r(\X_2; 0)=
3^{r-2}\Bigl(3 S_{r-2}(\X_2) - 2S_{1,r-3}(\X_2)\Bigr)\,.
\end{equation}
\end{proposition}
(In particular, $V_2(\X_2;0)=5$ and $V_3(\X_2;0)=9S_1(\X_2)$\,.)

\smallskip

\noindent
The proof of the proposition is given in the Appendix.

\bigskip

We now determine the specialization 
$F_r(\X_2\moins \D\moins\B_{r-2})$.

\begin{lemma}\label{div} The resultant $R(\X_2,\D+\B_{r-2})$ divides \
$F_r(\X_2\moins \D\moins \B_{r-2})$.
\end{lemma}
\proof
By \cite[Proposition 10]{P3}, we have
$$
F_r(x-\B_r)=R(x+\fbox{$2x$}+\fbox{$3x$},\B_r)\,,
$$
and making into $F_r(\X_2\moins \D\moins\B_{r-2})$
the substitutions: $x_1=0$ and $x_1=2x_2$, we get
$$
F_r(-\fbox{$2x_2$}-\B_{r-2})=R(0+0+0, \fbox{$2x_2$}
+\B_{r-2}+0)=0\,,
$$
and
$$
\aligned
F_r(x_2\moins\fbox{$2x_1$}\moins\fbox{$x_1\plus x_2$}\moins\B_{r-2})
=R(x_2\plus\fbox{$2x_2$}\plus\fbox{$3x_2$},
\fbox{$2x_1$}\plus \fbox{$x_1\plus x_2$}\plus\B_{r-2})\\
=R(x_2\plus \fbox{$2x_2$}\plus \fbox{$3x_2$},
\fbox{$2x_1$}\plus \fbox{$3x_2$}\plus \B_{r-2})=0\,. \ \ \ &
\endaligned
$$
Moreover, if $x_1\in \B_{r-2}$ and $\B_{r-3}:=\B_{r-2}\moins x_1$,
then $F_r(\X_2\moins \D\moins\B_{r-2})$ becomes
$$
\aligned
F_r(x_2\moins \fbox{$2x_1$} \moins \fbox{$2x_2$}
\moins &\fbox{$x_1\plus x_2$}\moins \B_{r-3})\\
&=R(x_2 \plus\fbox{$2x_2$}\plus\fbox{$3x_2$}, \fbox{$2x_1$}\plus
\fbox{$2x_2$}\plus \fbox{$x_1 \plus x_2$}\plus \B_{r-3})=0\,.
\endaligned
$$
These vanishings imply the assertion of the lemma.
\qed

\medskip

\noindent
We set
\begin{equation}
U_r(\X_2;\B_{r-2})=\frac{F_r(\X_2\moins \D\moins \B_{r-2})}{R(\X_2,\D+\B_{r-2})}\,.
\end{equation} 

Note that each variable
$b\in \B_{r-2}$ appears at most with degree $3$ in
$F_r(\X_2\moins \D\moins \B_{r-2})$, and hence at most
with degree $1$ in $U_r(\X_2;\B_{r-2})$. We have the following precise
recursive relation which follows from the observation that
the coefficient of $b_{r-2}^3$ in
$F_r(\X_2\moins \D\moins \B_{r-2})$
is equal to $F_{r-1}(\X_2\moins \D\moins \B_{r-3})$.

\begin{lemma}\label{Lur} For $r\ge 2$, we have
\begin{equation}
U_r(\X_2; \B_{r-2})= \sum_{i=0}^{r-2} \ U_{r-i}(\X_2; 0)
\ S_i(-\B_{r-2})\,.
\end{equation}
\end{lemma}

Let $\pi$ be the endomorphism of the $\bf C$-vector space of functions 
of $x_1,x_2$, defined by 
$$
\pi \bigl(f(x_1,x_2)\bigr) =
\frac{x_1 f(x_1,x_2)-x_2 f(x_2,x_1)}{x_1-x_2}\,.
$$
For any $i,j\in\N$, we have
\begin{equation}\label{pi}
\pi(x_1^jx_2^i)=S_{i,j}(\X_2)\,.
\end{equation}

\begin{proposition}\label{gr} The following identity holds for $r\ge 2$,
\begin{equation}\label{ilo} 
F_r(\X_2-\D)=-3^{r-2}R(X_2,\D) (x_1x_2)^{r-2}\bigl(3S_{r-2}(\X_2)-2S_{1,r-3}(\X_2)\bigr)\,.
\end{equation}
\end{proposition}
\proof
The identity is true for $r=2$.
To prove the assertion for $r\ge 3$, 
we compute in two different ways the action of $\pi$ on the multi-Schur function
(see \cite[1.4.7]{L} p.~9): 
\begin{equation}\label{multi}
S_{r,r;r}(\X_2+\fbox{$2x_1$}+\fbox{$3x_1$}-\D;x_1-\D)\,.
\end{equation}

Firstly, expanding (\ref{multi}), we have
$$
\aligned
&\pi \bigl(S_{r,r;r}(\X_2+\fbox{$2x_1$}+\fbox{$3x_1$}-\D;x_1-\D)\bigr)\cr
&=\pi \bigl(\sum_{j_1\le j_2 \le r} S_{j_1,j_2}(\fbox{$2x_1$}+\fbox{$3x_1$}) \ 
S_{r-j_2,r-j_1,r}(\X_2-\D;x_1-\D)\bigr)\cr
&=\pi \bigl(\sum_{j_1\le j_2 \le r} S_{j_1,j_2}(\fbox{$2$}+\fbox{$3$}) \ 
S_{r-j_2,r-j_1,r+j_1+j_2}(\X_2-\D;x_1-\D)\bigr)\cr
&=\sum_{j_1\le j_2 \le r} S_{j_1,j_2}(\fbox{$2$}+\fbox{$3$}) \ 
S_{r-j_2,r-j_1,r+j_1+j_2}(\X_2-\D)\cr
&=F_r(\X_2-\D)\,.
\endaligned
$$

\medskip

Secondly, we subtract $x_1$ from the arguments in the first two rows 
of (\ref{multi}) without changing the determinant (see \cite[Transformation Lemma 1.4.1]{L}):
\begin{align}\label{dis}
&S_{r,r;r}(\X_2+\fbox{$2x_1$}+\fbox{$3x_1$}-\D;x_1-\D)\cr
&=S_{r,r;r}(\X_2+\fbox{$3x_1$}-\fbox{$2x_2$}-\fbox{$x_1+x_2$};x_1-\D)\,.
\end{align}
Then the elements in the first two rows of the last column become zero, and we get 
the following factorization of the latter determinant in (\ref{dis}):
$$
S_{r,r}(x_2+\fbox{$3x_1$}-\fbox{$2x_2$}-\fbox{$x_1+x_2$}) \cdot S_r(x_1-\D)\,.
$$
Using the following two factorizations:
$$
S_{r,r}(x_2+\fbox{$3x_1$}-\fbox{$2x_2$}-\fbox{$x_1+x_2$})= 
-3^{r-2}(x_2-2x_1)(x_1x_2)^{r-1}(3x_1-2x_2)\,,
$$
and
$$
S_r(x_1-\D)=x_1^{r-2}x_2(x_1-2x_2)\,,
$$
we infer that
\begin{equation}\label{factor}
S_{r,r;r}(\X_2+\fbox{$2x_1$}+\fbox{$3x_1$}-\D;x_1-\D)=
-3^{r-2}R(\X_2,\D)(x_1x_2)^{r-2}x_1^{r-3}(3x_1-2x_2)\,.
\end{equation}
By (\ref{pi}), the result of applying $\pi$ to (\ref{factor}) is
$$
-3^{r-2}R(X_2,\D) (x_1x_2)^{r-2}\bigl(3S_{r-2}(\X_2)-2S_{1,r-3}(\X_2)\bigr)\,.
$$

Comparison of both these computations of $\pi$ applied to (\ref{multi}) 
yields the proposition.
\qed

\smallskip

In terms of $U_r$, we rewrite Proposition \ref{gr} into 

\begin{corollary}\label{Pu} \ For $r\ge 2$,
\begin{equation}\label{Ur0}
U_r(\X_2;0)=-3^{r-2}\bigl(3 S_{r-2}(\X_2)-2S_{1,r-3}(\X_2)\bigr)\,.
\end{equation}
\end{corollary}

\bigskip

Lemmas \ref{Lvr}, \ref{Lur}, Proposition \ref{Pv},
and Corollary \ref{Pu} imply Eq. (\ref{frhr}), and this 
finishes the proof of Theorem \ref{TA3}.

\bigskip

\section{Appendix: The Pascal starcaise}

We shall use the following variant of the Pascal triangle.
Consider an infinite matrix $P=[p_{s,t}]$ with rows and columns numbered
by $s,t=1,2, \ldots$. 

We assume that 
$p_{1,t}=p_{2,t}=0$ for $t\ge 2$, $p_{3,t}=p_{4,t}=0$
for $t\ge 3$, $p_{5,t}=p_{6,t}=0$ for $t\ge 4$ etc.
(Speaking less formally, $P$ is filled with 0's above the diagram
of the infinite partition  $(0,0,1,1,2,2,3,3,\ldots$) \ .)

The first column is an arbitrary sequence $v=(v_1,v_2,\ldots)$. In the case when 
this sequence is the sequence of coefficients of the Taylor expansion of a function $f(z)$, 
we write $P_f$ for the corresponding $P$.

To define the remaining $p_{s,t}$'s, we use the recursive formula
\begin{equation}
p_{s+1,t}= p_{s,t-1} + p_{s,t}.
\end{equation}
We visualize this definition by
$$ \begin{matrix}   a & b\\ & \square \end{matrix}   \qquad
\Rightarrow   \qquad \begin{matrix}   a & b\\ & a+b \end{matrix}
$$
We thus get the following {\it Pascal staircase} $P=[p_{i,j}]_{i,j=1,2,\ldots}$:
$$
\begin{array}{cccccc}
v_1 & 0 & 0 & 0 & 0 & \ldots \\
v_2 & 0 & 0 & 0 & 0 & \ldots \\
v_3 & v_2 & 0 & 0 & 0 & \ldots \\
v_4 & v_3\plus v_2 & 0 & 0 & 0 & \ldots \\
v_5 & v_4\plus v_3\plus v_2 & v_3\plus v_2 & 0 & 0 & \ldots \\
v_6 & v_5\plus v_4\plus v_3\plus v_2 & v_4\plus 2v_3\plus 2v_2 & 0 & 0 & \ldots \\
v_7 &  v_6\plus v_5\plus v_4\plus v_3\plus v_2 & v_5\plus 2v_4\plus 3v_3\plus 3v_2 & v_4\plus 2v_3\plus 2v_2& 0 & \ldots \\
\vdots & \vdots & \vdots & \vdots & \vdots &
\end{array}
$$

Given an integer $n\ge 0$, and an alphabet $\A$, we define the function
$W(n)=W(n,\A)$ by
\begin{equation}
W(n,\A)=\sum_{i,j} p_ {n+1-i,j+1}\, S_i(-\A)\, S_{j,n-i-j}(\X_2).
\end{equation}

\smallskip

The function $W(n,\A)$ is linear in the elements of the first 
column of $P$. Therefore it is sufficient to restrict to the case 
$v=(1,y,y^2,\ldots)$, i.e. to take \ $P= P_{1/(1-zy)}$ to determine it.

\begin{lemma} \ If $P= P_{1/(1-zy)}$ and $\A= \fbox{$x_1+x_2$}$, then $W(0)=1$ 
and for $n\ge 1$
\begin{equation}\label{lw}
W(n,\fbox{$x_1+x_2$}) = (y-1)y^{n-1} S_n(\X_2).
\end{equation}
\end{lemma}
\proof The entries contributing to $S_{k,n-k}(\X_2)$, 
where $k>0$ and $2k<n$ are, for some $a,b$, 
$$ 
\begin{matrix}  -a (x_1+x_2) S_{k-1,n-k}(\X_2) & -b(x_1+x_2)S_{k,n-k-1}(\X_2) \\
      &    (a+b) S_{k,n-k}(\X_2) 
      \end{matrix} 
      $$
and give $-a S_{k,n-k}(\X_2)-bS_{k,n-k}(\X_2) +(a+b)S_{k,n-k}(\X_2)=0$.

\smallskip

The entries contributing to $S_{k,k}(\X_2)$, where $k>0$ and $n=2k$ are, for some $a$,
$$ 
\begin{matrix}  -a (x_1+x_2) S_{k,k}(\X_2) & 0 \\
      &    a S_{k,k}(\X_2) 
      \end{matrix} 
      $$
and give $-a S_{k,k}(\X_2)+aS_{k,k}(\X_2)=0$.

\smallskip

Moreover, the first column contributes to $(y^n -y^{n-1})S_n(\X_2)$.
\qed

\bigskip

Taking now $\A=\fbox{$x_1+x_2$}+\B$ instead of $\fbox{$x_1+x_2$}$,
and using that 
$$ 
W(n,\A)=
\sum_{i,j,k} p_ {n+1-i-k,j+1}\, S_i\bigl(-\fbox{$x_1+x_2$}\bigr)\, 
S_{j,n-i-j-k}(\X_2) S_k(-\B) 
$$
$$ 
= \sum_k W\left(n-k,\fbox{$x_1+x_2$}\right) S_k(-\B)  
$$
$$ 
= (1-y^{-1}) \sum_k y^{n-k} S_{n-k}(\X_2) S_k(-\B) =  y^n ((1-y^{-1}) S_n( \X_2-y^{-1}\B)\,,  
$$
we get the following corollary.

\begin{corollary} \ For $P= P_{1/(1-zy)}$, $\B$ an arbitrary alphabet,
then (apart from initial values), we have
\begin{equation}
W(n,\fbox{$x_1+x_2$}+\B)=(y-1)y^{n-1} S_n(\X_2-y^{-1}\B)\,.
\end{equation}
\end{corollary}

We apply the corollary with $\B=\fbox{$2x_1$}+\fbox{$2x_2$}$. Expanding
\begin{multline*}
 S_n\left(\X_2 - y^{-1}(\fbox{$ 2x_1$}+\fbox{$ 2x_2$})\right)\\
= S_n(\X_2) - \frac{2x_1+2x_2}{y} S_{n-1}(\X_2)
 + 4\frac{x_1x_2}{y^2} S_{n-2}(\X_2)\,,
\end{multline*}
we get, for $n\ge 3$,
\begin{equation}
W(n,\D) = y^{n-2}(y-1)(y-2) S_{n}(\X_2))-2y^{n-3}(y-1)(y-2) S_{1,n-1}(\X_2)
\end{equation}
and initial conditions
$$
W(0) =1, \ \ W(1)= (y-3)S_1(\X_2), \ \ W(2)=(y-1)(y-2)S_2(\X_2)-2(y-3)S_{11}(\X_2)\,.
$$

\smallskip 

We come back to Proposition \ref{Pv}, and we take the Pascal staircase $P_f$ 
associated with the function
$$   
f=\frac{5-6z}{(1-z)(1-2z)(1-3z)}  =
- \frac{1/2}{1-z} - \frac{8}{1-2z}  + \frac{27/2}{1-3z}.
$$
Then for $P=P_f$, and $n=r-2$\,, the function $W(n,\D)$ is the function $V_r(\X_2;0)$.

\smallskip

We thus have to specialize $y$ into $1,2,3$ successively.
Apart from initial values, only $y=3$ contributes, 
and we get, for $n\ge 3$,
$$ 
W(n,\D)= 3^{n+1}S_n(\X_2)- 2\cdot 3^n S_{1,n-1}(\X_2)\,.
$$

This proves Proposition \ref{Pv}, checking the cases $r=2,3,4$ directly.

\bigskip

\noindent
{\bf Note}
As the referee of \cite{P3} points out, the Thom polynomials for Morin singularities
have been recently also studied -- using quite different methods --
by Feh\'er and Rim\'anyi in \cite{FR1}, and by B\'erczi and
Szenes in \cite{BSz}.


\bigskip

\begin{thebibliography}{99}\small
\addcontentsline{toc}{section}{\string\numberline{}References}

\bibitem{AVGL} V. Arnold, V. Vasilev, V. Goryunov, O. Lyashko:
\emph{Singularities. Local and global theory},
Enc. Math. Sci. vol. 6 (Dynamical Systems VI), Springer, 1993.

\bibitem{BFR} G. B\'erczi, L. Feh\'er, R. Rim\'anyi,
\emph{Expressions for resultants coming from the global theory of
singularities}, in: ``Topics in algebraic and noncommutative geometry'',
(L. McEwan et al. eds.), Contemporary Math. AMS {\bf 324} (2003), 63--69.

\bibitem{BSz} G. B\'erczi, A. Szenes,
\emph{Thom polynomials of Morin singularities},
arXiv: math.AT/0608285.

\bibitem{BR} A. Berele, A. Regev,
\emph{Hook Young diagrams with applications to combinatorics and to
representation theory of Lie superalgebras},
Adv. in Math. {\bf 64} (1987), 118--175.

\bibitem{D} J. Damon,
\emph{Thom polynomials for contact singularities},
Ph.D. Thesis, Harvard, 1972.

\bibitem{FK} L. Feh\'er, B. Komuves,
\emph{On second order Thom-Boardman singularities},
Fund. Math. {\bf 191} (2006), 249-264.

\bibitem{FR} L. Feh\'er, R. Rim\'anyi,
\emph{Calculation of Thom polynomials and other cohomological obstructions
for group actions}, in: ``Real and complex singularities (S\~ao Carlos
2002)'' (T. Gaffney and M. Ruas eds.), Contemporary Math. {\bf 354}
(2004), 69--93.

\bibitem{FR1} L. Feh\'er, R. Rim\'anyi,
\emph{On the structure of Thom polynomials of singularities},
Bull. London Math. Soc. {\bf 39} (2007), 541-549. 

%\bibitem{FP} W. Fulton, P. Pragacz,
%\emph{Schubert varieties and degeneracy loci},
%Springer LNM {\bf 1689} (1998).

%\bibitem{J} K. J\"anich,
%\emph{Symmetry properties of singularities of $C^{\infty}$-functions},
%Math. Ann. {\bf 238} (1979), 147--156.

\bibitem{Ka} M. \'E. Kazarian,
\emph{Characteristic classes of singularity theory},
in: ``The Arnold-Gelfand mathematical seminars: Geometry and singularity
theory'' (1997), 325--340.

\bibitem{Ka2} M. \'E. Kazarian,
\emph{Classifying spaces of singularities and Thom polynomials},
in: ``New developments in singularity theory'', NATO Sci. Ser. II Math.
Phys. Chem., {\bf 21}, Kluwer Acad. Publ., Dordrecht (2001), 117--134.

\bibitem{Kl} S. Kleiman,
\emph{The enumerative theory of singularities},
in: ``Real and complex singularities, Oslo 1976'' (P. Holm ed.) (1978),
297--396.

\bibitem{L} A. Lascoux,
\emph{Symmetric functions and combinatorial operators on polynomials},
CBMS/AMS Lectures Notes {\bf 99}, Providence (2003).

%\bibitem{L1} A. Lascoux,
%\emph{Addition of $\pm 1$: application to arithmetic},
%S\'eminaire Lotharingien de Combinatoire, {\bf B52a} (2004), 9 pp.

%\bibitem{LS} A. Lascoux, M-P. Sch\"utzenberger,
%\emph{Formulaire raisonn\'e de fonctions sym\'e\-triques},
%Universit\'e Paris 7, 1985.

\bibitem{M} I. G. Macdonald,
\emph{Symmetric functions and Hall-Littlewood polynomials},
Oxford Math. Monographs, Second Edition, 1995.

\bibitem{dPW} A. du Plessis, C. T. C. Wall,
\emph{The geometry of topological stability},
Oxford Math. Monographs, 1995.

\bibitem{P} P. Pragacz,
\emph{Enumerative geometry of degeneracy loci},
Ann. Sc. Ec. Norm. Sup. {\bf 21} (1988), 413--454.

\bibitem{P2} P. Pragacz,
\emph{Algebro-geometric applications of Schur $S$- and $Q$-polyno\-mials},
in: ``Topics in invariant theory'' -- S\'eminaire d'Alg\`ebre
Dubreil-Malliavin 1989-1990 (M-P. Malliavin ed.), Springer LNM {\bf 1478}
(1991), 130--191.

\bibitem{P23} P. Pragacz,
\emph{Thom polynomials and Schur functions I},
math.AG/0509234.

\bibitem{P4} P. Pragacz,
\emph{Thom polynomials and Schur functions: the singularities
$I_{2,2}(-)$},
Ann. Inst. Fourier {\bf 57} (2007), 1487--1508.

\bibitem{P3} P. Pragacz,
\emph{Thom polynomials and Schur functions: towards the singularities
$A_i(-)$,} 
in: ''Real and complex singularities - Sao Carlos 2006'', (M.~J. Saia and J. Seade eds.),
Contemporary Math. AMS {\bf 459}, (2008) 165-178.

\bibitem{PT} P. Pragacz, A. Thorup,
\emph{On a Jacobi-Trudi identity for supersymmetric polynomials},
Adv. in Math. {\bf 95} (1992), 8--17.

\bibitem{PW} P. Pragacz, A. Weber,
\emph{Positivity of Schur function expansions of Thom polynomials},
Fund. Math. {\bf 195} (2007), 85--95.

%\bibitem{PW2} P. Pragacz, A. Weber,
%\emph{Thom polynomials of invariant cones, Schur functions, and positivity},
%in: ``Algebraic cycles, sheaves, shtukas, and moduli'', (P. Pragacz ed.), Trends 
%in Mathematics, Birkh\"auser (2007) 117--129.

\bibitem{Rim2} R. Rim\'anyi,
\emph{Thom polynomials, symmetries and incidences of singularities},
Inv. Math. {\bf 143} (2001), 499--521.

\bibitem{RS} R. Rim\'anyi, A. Sz\"ucs,
\emph{Generalized Pontrjagin-Thom construction for maps with
singularities},
Topology {\bf 37} (1998), 1177--1191.

\bibitem{S} J. Stembridge,
\emph{A characterization of supersymmetric polynomials},
J. Algebra {\bf 95} (1985), 439-87.

\bibitem{T} R. Thom,
\emph{Les singularit\'es des applications diff\'erentiables},
Ann. Inst. Fourier {\bf 6} (1955--56), 43--87.

%\bibitem{W} C. T. C. Wall,
%\emph{A second note on symmetry of singularities},
%Bull. London Math. Soc. {\bf 12} (1980), 347--354.

\end{thebibliography}
\end{document}